\documentclass[12pt]{amsart}

\usepackage{amssymb}
\usepackage{upref}

\usepackage[lite]{amsrefs}

\newtheorem{thm}{Theorem}[section]

\theoremstyle{definition}

\numberwithin{equation}{section}

\newcommand{\secref}[1]{Section~\textup{\ref{#1}}}

\newcommand{\OO}{\mathcal O}

\renewcommand{\a}{\alpha}

\newcommand{\p}{\phi}
\newcommand{\x}{\xi}
\newcommand{\n}{\eta}
\renewcommand{\t}{\theta}
\newcommand{\Chi}{\raisebox{2pt}{\ensuremath{\chi}}}

\renewcommand{\P}{\Phi}

\DeclareMathOperator*{\spn}{span}
\DeclareMathOperator*{\clspn}{\overline{\spn}}

\renewcommand{\:}{\colon}
\newcommand{\<}{\langle}
\renewcommand{\>}{\rangle}

\begin{document}

\title[A theorem of Ionescu]{Bundles of $C^*$-correspondences
over directed graphs and a theorem of Ionescu}

\author{John Quigg}
\address{Department of Mathematics and Statistics
\\Arizona State University
\\Tempe, Arizona 85287}
\email{quigg@math.asu.edu}

\subjclass[2000]{Primary 46L08}

\keywords{directed graph, $C^*$-correspondence, graph
$C^*$-algebra, Cuntz-Pimsner algebra}

\date{\today}

\begin{abstract}
We give a short proof of a recent theorem of Ionescu which shows
that the Cuntz-Pimsner $C^*$-algebra of a certain correspondence
associated to a Mauldin-Williams graph is isomorphic to the graph
algebra.
\end{abstract}

\maketitle

\section{Introduction}

In recent decades the study of fractal geometry has led to the
introduction of \emph{graph-directed iterated function systems}
\cite{mauldin}, also known as \emph{Mauldin-Williams graphs}
\cites{edgar, io}. These are
finite
directed graphs of contractions
among compact metric spaces.
Recently, Ionescu \cite{io} associated a $C^*$-correspondence to
a given Mauldin-Williams graph, and proved that the resulting
Cuntz-Pimsner algebra is isomorphic to the graph $C^*$-algebra.
Ionescu constructs the isomorphism directly, extending the
contractions to the state spaces using Rieffels's theory of
Lipschitz metrics on state spaces \cite{rie:metric}.
Ionescu's result is perhaps surprising, and
illustrates important connections among
fractal geometry,
$C^*$-correspondences, and graph algebras.
Thus, we feel it will be useful to
show how Ionescu's theorem
can be quickly deduced from the elementary theory of graph
algebras.

We relax Ionescu's hypotheses somewhat:
whereas the directed graphs in \cite{io} are finite and have
neither sources nor sinks,
we only require the
graph to be row-finite with no sources. In fact, we only impose
these conditions to illustrate our method in its simplest form;
the general case could be handled with somewhat more effort.

Also, instead of contractions among metric spaces, we only
require continuous maps among locally compact Hausdorff spaces,
together with an equivariant surjection from the infinite path
space (see \secref{main} for details). For finite graphs,
such continuous maps together with an equivariant surjection
constitute a \emph{self-similarity structure} \cite{kig}.
As pointed out in \cites{kig, io},
every Mauldin-Williams graph gives rise to a self-similarity
structure.

\section{Preliminaries}

Let $E=(E^0,E^1,r,s)$ be a (directed) graph,
with vertices $E^0$, edges $E^1$, and range and source maps $r$
and $s$.
For $u,v\in E^0$ put $E^1_{uv}=\{e\in E^1\mid r(e)=u,s(e)=v\}$.
Warning:
we use the relatively new convention
(see \cites{katsura:class, rae:graph})
regarding the graph algebra $C^*(E)$: the generators go the same
way as the edges.
Thus, for example,
if
$e\in E^1_{uv}$ then $s_e^*s_e=p_v$ and $s_es_e^*\le p_u$, and
a finite path $e_1\cdots e_n$ in $E$ satisfies
$s(e_i)=r(e_{i+1})$ for $i=1,\dots,n-1$.

For simplicity we assume throughout that $E$ is row-finite and
has no sources,
meaning that each vertex receives a positive but finite number of
edges.

Let $E^*$ be the set of finite paths,
where vertices are regarded as paths of length $0$.
Let $E^\infty$ denote the set of infinite paths,
which under our hypotheses
is locally
compact Hausdorff when given the product topology.
Extend the source and range maps to
paths
in the obvious way, and
for $v\in E^0$
put $E^*_v=\{\a\in E^*\mid r(\a)=v\}$,
and similarly for $E^\infty_v$.

For $\a\in E^*$ let
$p_\a=s_\a s_\a^*$ be the range projection of the generator
$s_\a$.
Put $A_E=\clspn\{p_\a\mid \a\in E^*\}$, a commutative
$C^*$-subalgebra of $C^*(E)$.
It is folklore that
there is an isomorphism $\t\:C_0(E^\infty)\to A_E$
which takes the characteristic function of the set of
infinite paths starting
with a finite path $\a$
to
the generating projection $p_\a$.
We have $\t(C_0(E^\infty_v)=p_vA_E$.

Each $e\in E^1_{uv}$ gives rise to a continuous map
$\p^E_e\:E^\infty_v\to E^\infty_u$ via
$\p^E_e(\a)=e\a$.
For $f\in C_0(E^\infty_v)$ we have
\[
\t(f\circ \p^E_e)=s_e^*\t(f)s_e.
\]

\section{Ionescu's Theorem}
\label{main}

Suppose that for each $v\in E^0$ we have a $C^*$-algebra $A_v$,
and for each $e\in E^1_{uv}$ we have an $A_u-A_v$
correspondence $X_e$.
Let $A=\bigoplus_{v\in E^0}A_v$ be the $c_0$-direct sum.
Then each $X_e$ can be regarded as a correspondence over $A$;
let $X=\bigoplus_{e\in E^1}X_e$ be the direct sum of these
correspondences,
with operations
\[
(a\x)_e=a_{r(e)}\x_e,
\qquad
(\x a)_e=\x_ea_{s(e)},
\qquad
\<\x,\n\>_v=\sum_{s(e)=v}\<\x_e,\n_e\>
\]
for $a=(a_v)_{v\in E^0}\in A$ and $\x=(\x_e)_{e\in
E^1},\n=(\n_e)_{e\in
E^1}\in X$.
Then $X$ is a correspondence over $A$.
Let $\OO_X$ be the associated Cuntz-Pimsner algebra.

Here we are interested in correspondences arising as follows:
for each $v\in E^0$ let $T_v$ be a locally compact Hausdorff
space,
and put $A_v=C_0(T_v)$.
For each $e\in E^1_{uv}$ let $\p_e\:T_v\to T_u$ be a
continuous map.
and let $\p_e^*\:A_u\to A_v$
be the associated homomorphism, so that $A_v$ becomes a
$A_u-A_v$ correspondence $X_e$ with (right) Hilbert
module structure coming from the operations of $A_v$
and left module multiplication defined using $\p_e^*$.
Then $A=\bigoplus_{v\in E^0}A_v$ can be identified with $C_(T)$,
where $T$ is the disjoint union of $\{T_v\mid v\in E^0\}$.
Let $X=\bigoplus_{e\in E^1}X_e$ as above.

If the graph $E$ is finite, each $T_v$ is a compact metric
space, and each $\p_e$ is a contraction, we have a
Mauldin-Williams graph,
and $\OO_X$ as above was introduced in \cites{io,
bartholdi}.
Ionescu's proves \cite{io}*{Theorem 2.3} that $\OO_X$ is
isomorphic to the graph algebra $C^*(E)$.

Suppose we have a continuous map $\P\:E^\infty \to T$
satisfying $\P\circ \p^E_e=\p_e\circ \P$ for all $e\in E^1$.
If $E$ is finite and $\P$ is surjective, we have a
self-similarity structure \cite{kig}.

Ionescu's theorem follows from the following result:

\begin{thm}
[\cite{io}*{Theorem 2.3}]
With the above notation, suppose that
the continuous map
$\P\:E^\infty\to T$ is surjective. Then the Cuntz-Pimsner algebra
$\OO_X$ is isomorphic to the graph algebra $C^*(E)$.
\end{thm}

\begin{proof}
Define $\pi\:C_0(T)\to C^*(E)$ by $\pi(f)=\t(f\circ \P)$,
and for each $e\in E^1$ define $\psi_e\:X_e\to C^*(E)$ by
$\psi_e(\x)=s_e\pi(\x)$.
Routine calculations, using $p_v=\t(\Chi_{T_V}\circ \P)$ and
commutativity of $A_E$, show that the pair $(\psi_e,\pi)$ is a
covariant representation of the correspondence
$X_e$. We can form the direct sum $\psi=\bigoplus_{e\in
E^1}\psi_e$, giving a
Cuntz-Krieger Toeplitz representation
$(\psi,\pi)$
of the correspondence
$X$ in $C^*(E)$.
Since the range of $(\psi,\pi)$ contains the generators of
$C^*(E)$, an application of the Gauge-Invariant Uniqueness
Theorem \cite{katsura:corresp}*{Theorem 6.4}
shows that the associated homomorphism $\psi\times\pi\:\OO_X\to
C^*(E)$ is an isomorphism.
\end{proof}

\end{document}